\documentclass[10pt]{amsart}
\usepackage{verbatim}
\usepackage{amsmath,amsfonts,amscd,amssymb,epsfig,euscript,bm}
\usepackage{amsxtra,mathrsfs}
\usepackage{pdfsync}
\newtheorem{theorem}{Theorem}[section]

\theoremstyle{definition}
{}

\theoremstyle{remark}

\newcommand{\pbp}{\mathcal{\partial \bar{\partial}}}

\numberwithin{equation}{section}
\usepackage{hyperref}
\hypersetup{colorlinks,citecolor=blue,plainpages=false,hypertexnames=false}
\begin{document}
\title{A note on Demailly's approach towards a conjecture of Griffiths}
\author{Vamsi Pritham Pingali}
\address{Department of Mathematics, Indian Institute of Science, Bangalore 560012, India}
\email{vamsipingali@iisc.ac.in}
\begin{abstract} 
We prove that a ``cushioned" Hermitian-Einstein-type equation proposed by Demailly in an approach towards a conjecture of Griffiths on the existence of a Griffiths positively curved metric on a Hartshorne ample vector bundle, has an essentially unique solution when the bundle is stable. This result indicates that the proposed approach must be modified in order to attack the aforementioned conjecture of Griffiths.
\end{abstract}
\maketitle
\section{Introduction}\label{sec:intro}
The notion of ampleness/positivity is paramount in algebraic geometry. For a holomorphic line bundle, there is only one notion of differentio-geometric positivity, i.e., there is a smooth Hermitian metric whose curvature form is a K\"ahler form. By the Kodaira embedding theorem, it coincides with algebro-geometric ampleness. A holomorphic vector bundle $E$ is said to be Hartshorne ample if $\mathcal{O}_{E^*}(1)$ is an ample line bundle over $\mathbb{P}(E^*)$. There is no unique differentio-geometric notion of positivity of curvature $\Theta$ of a smooth Hermitian metric $h$. There are several
competing inequivalent notions. The most natural of these notions are \emph{Griffiths positivity} ($\langle v, \sqrt{-1}\Theta v\rangle$ is a K\"ahler form for all $v\neq 0$), \emph{Nakano positivity} (the bilinear form defined by $\sqrt{-1}\Theta$ on $T^{1,0}M \otimes E$ is positive-definite), and \emph{dual-Nakano positivity} (the Hermitian holomorphic bundle $(E^*,h^*)$ is Nakano \emph{negative}). Nakano positivity and dual-Nakano positivity imply Griffiths positivity and all three of them imply Hartshorne ampleness. A famous conjecture of Griffiths \cite{Griffiths} asks whether Hartshorne ample vector bundles  admit Griffiths positively curved metrics. This conjecture is still open. However, a considerable amount of work has been done to provide evidence in its favour \cite{Bo, Camp, dem2, Liu,  MT, Naumann, Pinchern, Um}.\\
\indent Relatively recently, Demailly \cite{dem2} proposed a programme to prove the aformentioned conjecture of Griffiths for a holomorphic rank$-r$ vector bundle $E$ on a compact K\"ahler manifold $(X,\omega_0)$. In fact, if Demailly's method works, it will end up proving a stronger conjecture : Do Hartshorne-ample bundles admit dual-Nakano positively curved metrics ? Demailly's approach involves solving a family (depending on a parameter $0\leq t \leq 1$) of vector bundle Monge-Amp\`ere equations (distinct from the one introduced in \cite{Pin}) in conjunction with ``cushioned" Hermitian-Einstein-type equations (Theorem 2.17 in \cite{dem2}):
\begin{gather}
\det_{TX \otimes E^*} \left (\Theta_{h_t}+(1-t)\alpha \omega_0 \otimes I_{E^{*}} \right ) ^{1/r} = f_t \frac{(\det h_0)^{\mu}}{(\det h_t)^{\lambda}} \omega_0^n,  \\
\left (\sqrt{-1}F_{h_t} - \frac{\sqrt{-1}}{r} \mathrm{tr}F_{h_t}\right) \omega_0^{n-1} = -\epsilon\frac{(\det h_0)^{\mu}}{(\det h)^{\mu}} \ln \left (\frac{hh_0^{-1}}{\det(hh_0^{-1})^{1/r}} \right )\omega_0^n,
\label{eq:demeq}
\end{gather}
where $h_0$ is a smooth background Hermitian metric, $\mu, \lambda \geq 0$ are fixed constants, $\alpha>0$ is a large enough constant so that $\Theta_{h_0}+\alpha \omega$ is dual-Nakano positively curved, and $f_t>0$ are smooth positive functions. We focus on the cushioned Hermitian-Einstein-type equation in the following theorem.
\begin{theorem}
Let $E$ be an $\omega_0$-stable rank$-r$ holomorphic bundle on $X$. Let $H_0$ be a Hermitian-Einstein metric on $E$ with respect to $\omega_0$, that is, $\sqrt{-1}F_{H_0}\omega_0^{n-1}=\lambda \omega_0^n$. Let $h$ be a smooth metric on $E$ solving the following cushioned Hermitian-Einstein equation for given parameters $\epsilon \geq 0, \mu \geq 0$.
\begin{gather}
\left (\sqrt{-1}F_h - \frac{\sqrt{-1}}{r} \mathrm{tr}F_h\right) \omega_0^{n-1} = -\epsilon\frac{(\det H_0)^{\mu}}{(\det h)^{\mu}} \ln \left (\frac{hH_0^{-1}}{\det(hH_0^{-1})^{1/r}} \right )\omega_0^n,
\label{eq:viscousHE}
\end{gather}
where $h, H_0$ are matrices (any holomorphic trivialisation will do). Then $h=H_0e^{-f}$ for some smooth function $f$.
\label{thm:unique}
\end{theorem}
\indent As a result, if we consider the system of the vector bundle Monge-Amp\`ere equation and the cushioned Hermitian-Einstein-type equation on an $\omega_0$-stable ample $E$, and if solutions exist all the way till $t=1$, the final $t=1$ solution, by virtue of the fact that it satisfies the cushioned Hermitian-Einstein-type equation, has to be of the form $H_0 e^f$. This condition might be a strong restriction (which is unlikely to be met owing to \cite{Lubke,ST} without a restriction on the second Chern character). On the other hand, if we replace $\omega_0$ by say $(1-t)\omega_0 + t \sqrt{-1}tr(F_{h_t})$ (or the choice in Section 2.19 in \cite{dem2} for instance), the above argument will not be applicable and there might be some hope for the approach to yield an affirmative solution to the Griffiths Conjecture. 
\subsection*{Acknowledgements}
This work is partially supported by grant F.510/25/CAS-II/2018(SAP-I) from UGC (Govt. of India), and a MATRICS grant MTR/2020/000100 from SERB (Govt. of India). The author thanks Jean-Pierre Demailly for fruitful discussions and for encouraging me to write up this note.
\section{Proof of uniqueness}\label{sec:proof}
\indent In a holomorphic trivialisation, our conventions are : $\langle v, w \rangle_H = v^T H \bar{w}$, if $g$ is an endomorphism then $g.s = [g]^T\vec{s}$, $\nabla s = ds +A^Ts$, $A= \partial H H^{-1}$, $F= dA-A\wedge A=\bar{\partial} A$, and $\nabla g = dg + [g,A]$. \\
\indent The proof is motivated by a similar one by Donaldson for Riemann surfaces. In general, $h=qH_0$ where $q$ is some smooth $H_0$-Hermitian positive-definite endomorphism of $E$. We decompose $q$ further as $q=e^{-f} g$ where $\det(g)=1$ and $f$ is a smooth function. Thus, $F_h = F_{H_0} + \pbp f + \bar{\partial}(\partial_0 g g^{-1})$. The trace-free part of the curvature is $F_h^{\circ} = F_0 ^{\circ}+\bar{\partial}(\partial_0 g g^{-1})$. Substituting these expressions in \ref{eq:viscousHE} and using the fact that $H_0$ is Hermitian-Einstein with respect to $\omega_0$, we get 
\begin{gather}
\sqrt{-1}\bar{\partial}(\partial_0 g g^{-1}) \omega_0^{n-1} = -\epsilon e^{r\mu f} \ln g \omega_0^n.
\label{eq:gviscousHE}
\end{gather}
Now we compute
\begin{gather}
\frac{1}{2}\sqrt{-1} \bar{\partial} \partial tr(g^2)\omega_0^{n-1} = \sqrt{-1} \bar{\partial} tr(g\partial_0 g) \omega_0^{n-1} = \sqrt{-1} tr(\bar{\partial} g \partial_0 g)\omega_0^{n-1} +\sqrt{-1} tr(g\bar{\partial} \partial_0 g)\omega_0^{n-1} \nonumber \\
=\sqrt{-1} tr(\bar{\partial} g \partial_0 g)\omega_0^{n-1}-\epsilon e^{r\mu f} tr(g^2 ln g) \omega_0^n - \sqrt{-1}tr(g\partial_0 g g^{-1} \bar{\partial} g) \omega_0^{n-1}\nonumber \\
\leq -\epsilon e^{r\mu f} tr(g^2 ln g) \omega_0^n.
\end{gather}
Note that $tr(g^2 \ln g) = \sum_i \lambda_i ^2 \ln (\lambda_i)$ where $\lambda_i>0$ are the eigenvalues of $g$ such that $\lambda_1 \leq \lambda_2\ldots$. The product of the $\lambda_i$ is $1$. Thus, $$\displaystyle \sum_{i \ \vert  \ \lambda_i <1} \vert \ln (\lambda_i) \vert = \sum_{i \ \vert \ \lambda_i >1} \ln (\lambda_i)$$ which implies that 
\begin{gather}
\displaystyle \sum_{1\leq i \leq p \ \vert  \ \lambda_p \leq 1, \ \lambda_{p+1} >1}  \lambda_i^2 \vert \ln (\lambda_i) \vert \leq \sum_{1\leq i \leq p \ \vert  \ \lambda_p \leq 1, \ \lambda_{p+1}>1}  \lambda_p^2 \vert \ln (\lambda_i) \vert = \lambda_p^2\sum_{i =p+1}^n \ln (\lambda_i) \nonumber \\
 \leq \sum_{i=p+1}^n \lambda_i ^2 \ln (\lambda_i).
\end{gather}
 Therefore, $tr(g^2 \ln g) \geq 0$ and hence 
\begin{gather}
\frac{1}{2}\sqrt{-1} \bar{\partial} \partial tr(g^2)\omega_0^{n-1} \leq 0.
\end{gather}
The strong maximum principle then implies that actually
\begin{gather}
\frac{1}{2}\sqrt{-1} \bar{\partial} \partial tr(g^2)\omega_0^{n-1} =0,
\end{gather}
and hence $g=I$. \qed

\end{document}